\documentclass[a4paper]{article}
\usepackage{latexsym}
\usepackage{amssymb}
\usepackage{authblk}

\usepackage{xcolor} 
\newtheorem{Theorem} {Theorem} [section]
\newtheorem{Proposition} [Theorem] {Proposition}
\newtheorem{Corollary} [Theorem] {Corollary}
\newcommand{\Proof}{ \noindent{\bf Proof.}\quad }
\newcommand{\qed}{\hfill$\Box$\medskip}
\newcommand{\Ff}{{\mathbb F}}
\newcommand{\1}{{\bf 1}}
\newcommand{\eps}{\varepsilon}

\newcommand{\chisub}{\chi^{}}
\newcommand{\divides}{\,|\,}
\newcommand{\notdivides}{\nmid}
\newcommand{\maysplit}{\allowbreak}

\title{Cliques in Paley graphs of square order and in Peisert graphs}
\author[2]{Andries E. Brouwer}
\author[1]{Sergey Goryainov}
\author[3]{Leonid Shalaginov}
\author[4]{Chi Hoi Yip}

\affil[1]{School of Mathematical Sciences, Hebei International Joint Research Center for Mathematics and Interdisciplinary Science, Hebei Key Laboratory of Computational Mathematics and Applications, Hebei Workstation for Foreign
Academicians, Hebei Normal University,
P. R. China}
\affil[2]{Brouwer-Mortensen Institute for Retired Mathematicians and Artists,
Amsterdam}
\affil[3]{Chelyabinsk State University,
Russia}
\affil[4]{School of Mathematics, Georgia Institute of Technology,
USA}

\affil[${}$]{
\vskip -0.2cm
{\tt aeb@cwi.nl},
{\tt sergey.goryainov3@gmail.com}, 
{\tt 44sh@mail.ru},
{\tt cyip30@gatech.edu}}

\date{\today}

\begin{document}
\maketitle

\begin{abstract}
We study maximal cliques in the collinearity graphs of Desarguesian nets,
give some structural results and some numerical information.
In particular, we show for Desarguesian nets that the set consisting of
a point $x$ together with all its neighbors on a line $L$ (with $x$ not on $L$)
is contained in a unique maximal clique $C_{x,L}$ and determine the sizes
and automorphism groups of such maximal cliques $C_{x,L}$ in all cases.
\end{abstract}

\section{Introduction}
Throughout the paper, let $p$ be a prime and $q,r$ be prime powers. We use $\mathbb{F}_q$ to denote the finite field with $q$ elements.

Let $F$ be a finite field and $D \subseteq F \cup \{\infty\}$ a set
of `directions'.
Let $\Gamma(F,D)$ be the graph with vertex set $F \times F$, where
two vertices are adjacent when the line joining them has direction
in $D$. The Paley graphs of square order are examples of such graphs.
We study the cliques in these graphs $\Gamma(F,D)$.

\medskip
More precisely, in the following subsections of the present section,
we define Paley and Peisert graphs and survey what is known
about their maximal and second maximal cliques, and automorphism groups.
In Section \ref{cl-in-nets} we define the maximal cliques $C_{x,L}$
of Desarguesian nets, and find their sizes and automorphism groups.
In Section \ref{numerical} we give numerical results.

\medskip
Our graphs are finite and undirected, without loops.

For a graph $\Gamma$ with vertex $x$ we denote the set of neighbors
of $x$ by $\Gamma(x)$, and we write $x^\perp := \{x\} \cup \Gamma(x)$,
and $A^\perp := \bigcap_{a \in A} a^\perp$ for a set of vertices $A$.

For definition and properties of strongly regular graphs,
see \cite{BvM}.

\subsection{Paley graphs}
Let $q = 4t+1$ be a prime power.
The {\it Paley graph} $P(q)$ of order $q$ is the graph with vertex set
$\Ff_q$, where two field elements
are adjacent when their difference is a square in the field.
This graph is undirected, since $-1$ is a square in $\Ff_q$.
This graph is strongly regular with parameters
$(v,k,\lambda,\mu) = (4t+1,2t,t-1,t)$ and spectrum
\(\{ [2t]^1,\ \left[\frac{-1\pm\sqrt{q}}{2}\right]^{2t}\}\). It is self-complementary (isomorphic to its complement).

\subsection{Peisert graphs}
The {\it Peisert graph} $P^*(q)$ of order $q = p^{2e}$ where
$p \equiv 3$ (mod 4) is the graph with vertex set $\Ff_q$,
where two field elements are adjacent when their difference is $\beta^i$
with $i \equiv 0,1$ (mod 4), where $\beta$ is a fixed primitive element
in $\Ff_q$.
This graph is undirected, since $-1$ is a 4-th power in $\Ff_q$.
This graph is strongly regular with the same parameters as the
Paley graph of order $q$.
Peisert \cite{Peisert01} showed that a self-complementary graph
with automorphism group that is transitive on vertices and edges is
either a Paley graph, or a Peisert graph, or is isomorphic to a certain
sporadic graph $P^{**}(23^2)$ on $23^2$ vertices.
(For this last graph, see \cite[\S10.70]{BvM}.)

\medskip
{\footnotesize
In the literature some graphs are called `of Peisert type'
where Peisert graphs are not always `of Peisert type' and graphs
`of Peisert type' are not always self-complementary. A strange term,
better avoided. \par}

\subsection{Nets}
A {\it net} of degree $m$ and order $n$ is a
partial linear space on $n^2$ points with $mn$ lines
of size $n$ partitioned into $m$ parallel classes.
Clearly $0 \le m \le n+1$.

If $m = n+1$, then the net is an affine plane $AG(2,n)$
(that is, a 2-$(n^2,n,1)$ design), and the collinearity graph is a
complete graph.

If $0 < m < n+1$, then the collinearity graph is strongly regular
with parameters
$(v,k,\lambda,\mu) = (n^2,m(n-1),\maysplit (m-1)(m-2)+n-2,m(m-1))$.

If $m = (n+1)/2$, then the graph has the same parameters as the
complementary graph. (Such strongly regular graphs are known as
conference graphs.)

The eigenvalues of a strongly regular net graph are
$m(n-1)$, $n-m$ and $-m$ with multiplicities
$1$, $m(n-1)$ and $(n+1-m)(n-1)$, respectively.
A basis for the $(n-m)$-eigenspace is given by the vectors
$n \chisub_L-\1$ (where $\chisub_L$ is the characteristic function
of the line $L$, and $\1$ the all-1 vector), if we pick $n-1$ lines
from each parallel class.

If $n = r$ is a prime power, and $0 \le m \le n+1$, then we
can construct a net of degree $m$ and order $n$ by taking
$m$ parallel classes in the Desarguesian affine plane of order $r$. We shall refer to such a net as a Desarguesian net.
The Paley graphs of order $r^2$, and the Peisert graphs
of order $r^2$ with $r \equiv 3$ (mod 4)
are special cases of this construction (with $m = (n+1)/2$).

\subsection{Automorphism groups}
Let $q = p^e$, where $p$ is prime. The full automorphism group
of the Paley graph $P(q)$ of order $q$ consists of the maps
$x \mapsto ax^\sigma+b$, where $a,b \in \Ff_q$, $a$ a nonzero square,
and $\sigma = p^i$ with $0 \le i < e$ (Carlitz \cite{Carlitz60}).
It has order $eq(q-1)/2$.

The subgraph $\Pi$ of $P(q)$ induced on the neighbors of 0 has full
automorphism group consisting of the maps $x \mapsto a x^{\pm\sigma}$
with $a,\sigma$ as before (Muzychuk \& Kov\'acs \cite{MuzychukKovacs05}).
It has order $e(q-1)$ for $q > 9$ and is half as large for $q=5,9$.

\medskip
Let $q = p^f$, where $p \equiv 3$ (mod 4) is prime, and $f$ is even.
The full automorphism group of the Peisert graph $P^*(q)$ of order $q$
has order $fq(q-1)/4$ for $q \ne 3^2,7^2,3^4$, and is $2,3,6$
times as large in these three cases (Peisert \cite{Peisert01}).

\subsection{Maximum Cliques}
A {\it clique} is a complete subgraph of a graph.
For strongly regular graphs a standard upper bound for the size
of cliques is the Delsarte-Hoffman bound (cf.~\cite{BvM}).
For a strongly regular net graph this bound is $n$ (independent of $m$).
Of course this bound is attained by the lines. In particular,
for Paley graphs of order $r^2$, and for Peisert graphs of order $q = r^2$
with $r \equiv 3$ (mod 4), the subfield $\Ff_r$ is a clique of maximum size.
There remains the question about the maximum clique size for Peisert graphs
of order $q = s^4$. For $q = 3^4$ and $q = 7^4$ the answers turn out to be
9 and 17, respectively (Mullin \cite{Mullin09}). 

Blokhuis \cite{Blokhuis84} showed for Paley graphs of square order $q = r^2$
that the only cliques of size $r$ are lines of the underlying affine plane
of order $r$ (images $aF+b$ of the subfield $F = \Ff_r$, where
$a$ is a nonzero square). That is, $\Ff_r$ is the only subset of size $r$ of
$\Ff_{r^2}$ containing 0, 1 in which any two elements differ by a square.
More generally, Sziklai \cite{Sziklai99} showed for $r$ a prime power
and $d \divides (r+1)$ that $\Ff_r$ is the only subset of size~$r$~of
$\Ff_{r^2}$ containing $0, 1$ in which any two elements differ
by a $d$-th power. We also refer to a recent paper of Yip \cite{Y26} for an extension of these two theorems.

Such a classification is not known yet for Peisert graphs.
Asgarli \& Yip \cite{AsgarliYip22} proved some partial results.
They also showed that if $s > 3$ then $\Ff_s$ is a maximal clique
in $P^*(s^4)$.

In \cite{BvM}, \S8.4.2 there is some information on cliques in
general net graphs.
In particular it is observed that when $n > (m-1)^2$ each edge lies in a
unique $n$-clique, so that there are no other $n$-cliques than the lines
of the net.
On the other hand, if $n = (m-1)^2$ is the square of a prime power,
then affine Baer subplanes are examples of $n$-cliques that are not lines.

\subsection{Second Largest Maximal Cliques}
After the largest cliques, the second largest and the smallest maximal cliques
are of interest. In any strongly regular graph, if $C$ is a clique with
a size meeting the Delsarte-Hoffman bound, every point $x$ outside $C$
has the same number of neighbors inside.
In the case of a net graph this number is $m-1$. In particular, for Paley or
Peisert graphs of order $q = r^2$, each point $x$ outside an $r$-clique $C$
has $\frac{r-1}{2}$ neighbors inside. Now $\{x\} \cup (x^\perp \cap C)$
is an $\frac{r+1}{2}$-clique. Baker et al.~\cite{Baker-et-al96} showed that
these cliques are maximal in $P(r^2)$ when $r \equiv 1$ (mod 4), and that
$\{x,x^r\} \cup (x^\perp \cap C)$ is maximal in $P(r^2)$ when
$r \equiv 3$ (mod 4) and $C = \Ff_r$.
(Below we give an alternative proof.)

%
One conjectures that this size ($\frac{r+1}{2}$ if $r \equiv 1$ (mod 4)
or $\frac{r+3}{2}$ if $r \equiv 3$ (mod 4)) is the second largest size
for cliques in $P(r^2)$, and that there are only two orbits of cliques
of this size when $r \ge 25$. This conjecture has been checked
for $r \le 109$.

These two orbits arise as follows. One orbit consists of the Baker et al.\ cliques
\(\{x\}\cup (x^\perp \cap C)\) (or \(
\{x,x^q\}\cup (x^\perp \cap C) \) when \(q\equiv 3 \pmod{4}\)). As observed above, if \(C\) is a clique in $P(q)$ containing \(0\), then applying the automorphism $z \mapsto z^{-1}$ on the subgraph \(\Pi\) yields another clique of the same size, namely
\(\{0\}\cup \{c^{-1}: c\in C,\ c\neq 0\}.
\)
Propositions~\ref{Goryainov} and~\ref{Goryainov} below give a more symmetric description of this second orbit.

\medskip
There is a more symmetric description of these latter cliques.

\begin{Proposition} \label{Goryainov}
{\rm (Goryainov et al.~\cite{Goryainov-et-al18})}
Let $\beta$ be primitive in $\Ff_{r^2}$, and put $\omega := \beta^{r-1}$.
Let $Q := \langle \omega^2 \rangle$.
If $r \equiv 1~({\rm mod}~4)$ the set $Q$ is a maximal coclique
of size $(r+1)/2$ in $P(r^2)$.
If $r \equiv 3~({\rm mod}~4)$ the set $Q \cup \{0\}$ is a maximal
clique of size $(r+3)/2$ in $P(r^2)$.
\end{Proposition}
\Proof
Put $\eps = \beta^{(r+1)/2}$, so that $\eps^r = -\eps$.
Let $N \,:\, \Ff_{r^2} \to \Ff_r$ be the norm.
Then $N(x+y\eps) = x^2-dy^2$ where $d = \eps^2 \in \Ff_r$.
We see that $\langle \omega \rangle$ is the set of points
on the conic $x^2-dy^2 = 1$, and $Q$ consists of half of
the points on this conic.
Let $\omega^i = x+y\eps$ with $x,y \in \Ff_r$. Then
$N(\omega^{2i}-1) = ((x+y\eps)^2-1)((x-y\eps)^2-1) = -4dy^2$.
Since $-d$ is a square in $\Ff_r$ if and only if $r \equiv 3~({\rm mod}~4)$,
the given sets are (co)cliques as claimed. Maximality follows
from the following proposition.~\qed

\begin{Proposition} \label{Goryainov2}
{\rm (Goryainov et al.~\cite{Goryainov-et-al22})}
Let $L$ be the clique $\Ff_r$ in $P(r^2)$.
The map $z \mapsto \eps^{-1} (1 + \frac{2}{z-1})$,
$1 \mapsto \eps^{-1}$ maps $Q$ (resp.~$Q \cup \{0\}$) onto
$\{x\} \cup (x^\perp \cap L)$ (resp.~$\{x,x^r\} \cup (x^\perp \cap L)$),
where $x = \eps^{-1}$.
\end{Proposition}
\Proof
The given map maps 0, 1, $\eta = x+y\eps \in \langle \omega \rangle$ to
$-\eps^{-1}$, $\eps^{-1}$, and $\frac{y}{x-1}$, respectively.
Let $C$ be a maximal (co)clique containing $Q$ (resp.~$Q \cup \{0\}$).
Then $1 \in C$, so $z \mapsto \frac{1}{z-1}$ and
$z \mapsto 1 + \frac{2}{z-1}$ preserve adjacency on $C$,
so $z \mapsto \eps^{-1} (1 + \frac{2}{z-1})$
flips or preserves adjacency when $r \equiv 1~({\rm mod}~4)$ or
$r \equiv 3~({\rm mod}~4)$.~\qed

\medskip
Let $r = p^e$.
For $r \ge 9$, the Baker et al.~cliques have stabilizers
(in ${\rm Aut}\,P(r^2)$) of order $2e$ if $r \equiv 1$~(mod 4),
and $4e$ if $r \equiv 3$ (mod 4).
For $r \ge 5$, the Goryainov et al.~cliques have stabilizers
of order $e(r+1)$.

\section{Maximal cliques in Desarguesian nets} \label{cl-in-nets}
Throughout the section, assume that $p$ is the characteristic of the field
we are working on. Let a Desarguesian net of order $r$ and degree $m$ be
a net that is a subnet of the Desarguesian affine plane of order $r$
obtained by choosing $m$ directions. Its collinearity graph can be taken
to have vertex set $\Ff_{r^2}$, where two vertices are adjacent when
their difference lies in the set $S$, a union of $m$ cosets of $\Ff_r^*$
in $\Ff_{r^2}^*$.

We will repeatedly use the following standard facts about the Desarguesian affine plane $AG(2,r)$: translations and homotheties preserve directions and hence induce automorphisms of the corresponding net; lines in the same parallel class are translates of one another; and subgroups of ${\rm AGL}(1,r)$ act on a line in the usual affine manner. For background on these facts, see, for example, Hughes--Piper \cite{HughesPiper73} and Brouwer--Van Maldeghem \cite{BvM}.

Any clique in this graph is a set determining at most $m$ directions.
Conversely, an arbitrary subset of $\Ff_{r^2}$ that determines
$m$ directions is a clique in the net defined by these directions.
Blokhuis et al.~\cite{Blokhuis-et-al99} and Ball \cite{Ball03}
determined the possibilities for the number of directions $m$ determined
by a maximum clique (of size $r$). For $m \le \frac{r+1}{2}$ they show
that the clique must be a subspace.
(For definitions and details, see loc. cit.)
Applying this, Goryainov \& Yip \cite{GoryainovYip24} determined
the smallest possible degree $m$ for which there are maximum cliques
other than lines. They proved that if $r$ is prime, then $m = \frac{r+3}{2}$; if $r = s^e$ with $s$ maximal such that $e > 1$, then $m = s^{e-1}+1$ \cite[Theorem 1.3]{GoryainovYip24}. They also showed that for $e \le 3$, there is a unique example \cite[Theorem 1.5]{GoryainovYip24}.

\subsubsection*{Cliques $\{x\} \cup (x^\perp \cap L)$
are contained in unique maximal cliques $C_{x,L}$}
We generalize the Baker et al.~result, and show that if $L$ is a line
of the net, and $x$ a point outside, then $\{x\} \cup (x^\perp \cap L)$
is contained in a unique maximal clique.

\begin{Theorem} \label{netcliq}
Let $r$ be a power of the prime $p$, and consider
a Desarguesian net of order $r$ and degree $m$.
If $L$ is a line and $x$ a point outside, then the clique
$\{x\} \cup (x^\perp \cap L)$ in the collinearity graph of the net
is contained in a unique maximal clique $C_{x,L}$. Moreover, when
$p \nmid (m-1)$, if $\{x\} \cup (x^\perp \cap L)$ is not maximal,
then $C_{x,L}$ is contained in the union $L \cup M$ of two lines,
where $M$ is a unique line passing through $x$.
\end{Theorem}

\Proof
Let $A = x^\perp \cap L$. We prove that $C_{x,L} = (\{x\} \cup A)^\perp$.
It suffices to show that this is a clique, that is, that any two points
$v,w$ of $(\{x\} \cup A)^\perp$ are adjacent.

If $v$ or $w$ lies in $\{x\} \cup A$, there is nothing to prove.
Otherwise, $x,v,w$ are distinct points outside $L$.
We show that $vw$ is parallel to a net line on~$x$.

\smallskip
If the line $vw$ is parallel to $L$, it is a net line, and $v,w$
are adjacent. 

\smallskip
Let $T_{a,s} \colon z \mapsto a+s(z-a)$ be the homothety with center $a$
and dilatation factor $s$, where $s \in \Ff_r^*$. This map preserves
the directions of lines, and hence is an automorphism of the net.
This map also preserves each line on $a$.
If $a$ lies on $L$ and $y = T_{a,s}(x)$, $A = x^\perp \cap L$,
$B = y^\perp \cap L$, then $T_{a,s}$ maps $A$ onto $B$. In particular,
if $y \in (\{x\} \cup A)^\perp$ and $y\notin L$, then $A$ is invariant
under $T_{a,s}$.

\smallskip
Suppose that \(xv\) is parallel to \(L\) and \(xw\), \(vw\) are not. Now
\[
A=x^\perp\cap L=v^\perp\cap L=v-x+(x^\perp\cap L)=v-x+A,
\]
so \(A\) is invariant under the translation by \(v-x\). Let the affine lines \(xw\) and \(vw\)
meet \(L\) in points \(b\) and \(c\), respectively. Since \(w^\perp\cap L=A\), we see that \(A\)
is invariant under the homothety \(T_{b,t}\), where \(w=T_{b,t}(x)\) and \(t\neq 1\). Put
$d=b+(v-x)\in A.$
 To compute \(c\), choose affine coordinates so that \(L=\{(u,0):u\in\Ff_r\}\) and \(x=(0,1)\),
identifying points of \(L\) with their first coordinates. Then
$v=(d-b,1)$ and $w=((1-t)b,t).$ The line through \(v\) and \(w\) meets \(L\) when its second coordinate is \(0\), and this gives $c=b+\frac{t}{t-1}(d-b).$
Hence the line through \(x\) parallel to \(vw\) meets \(L\) in
$c-(v-x)=b+\frac{1}{t-1}(d-b).$
Now apply the affine relabeling \(z\mapsto (z-b)/(d-b)\) on \(L\). Since \(b\neq d\), this is
well defined, preserves the relevant invariance properties of \(A\), and sends \(b,d\) to \(0,1\),
respectively. Thus we may assume that \(b=0\) and \(d=1\). Now \(A\) contains \(0,1\), is
invariant under \(z\mapsto tz\) and under \(z\mapsto z+1\), and hence contains every value
\(p(t)\) for polynomials \(p\) with coefficients in \(\Ff_p\), and therefore the field
\(\Ff_p(t)\). It follows that
$\frac{1}{t-1}\in A,$ as desired.

\smallskip

The case where $xw$ is parallel to $L$, and $xv,vw$ are not,
follows similarly.

\smallskip
Finally suppose none of the lines $xv, xw, vw$ is parallel to $L$,
so that they meet $L$ in distinct points $a, b, c$, respectively.
Now $A$ is invariant under the homotheties $T_{a,s}$ and $T_{b,t}$,
where $v = T_{a,s}(x)$ and $w = T_{b,t}(x)$ and $s \ne t$.
Let $G = \langle T_{a,s},T_{b,t} \rangle$. Then $G$ contains
the translation $z \mapsto z+b-a$.

\smallskip
{\footnotesize
This is well-known in the theory of affine planes,
or that of Frobenius groups. See, e.g., \cite{HughesPiper73},
Theorem 4.25 and Corollary 1.
\par}

Computation shows that $c = a+\smash{\frac{s(1-t)}{s-t}(b-a)}$,
so that the line through $x$ parallel to $vw$ hits $L$
in $T_{a,s}^{-1} (c) = a + \frac{1-t}{s-t}(b-a)$.
After appropriate relabeling we may take $a=0$, $b=1$.
Now $A$ contains $0,1$, is invariant under $z \mapsto z+1$
and under $z \mapsto sz$ and under $z \mapsto 1+t(z-1)$, hence
under $z \mapsto tz$, and therefore $A$ contains the field $\Ff_p(s,t)$.

\smallskip
{\footnotesize
%
%
(Note that $\Ff_p(s,t) = \Ff_p(u)$ if $s = \alpha^i$, $t = \alpha^j$,
$u = \alpha^h$, where $\alpha$ is primitive in $\Ff_r^*$ and
$h = {\rm gcd}(i,j)$.)
\par}

\smallskip
It follows that $a+\frac{1-t}{s-t}(b-a)=\frac{1-t}{s-t} \in A$, as desired.

\medskip

For the final assertion, assume that \(p\nmid(m-1)\) and that \(\{x\}\cup A\) is not maximal.
Then there exists \(y\in C_{x,L}\setminus (\{x\}\cup A)\). As in the proof above, the line \(xy\)
cannot be parallel to \(L\), since otherwise \(A\) would be invariant under a nontrivial
translation, forcing \(p\mid (m-1)\). Hence \(xy\) meets \(L\) in a point \(b\), and \(A\) is
invariant under a homothety \(T_{b,t}\) with \(t\neq 1\). Summing over \(A\), we obtain
$\sum_{a\in A}(a-b)=t\sum_{a\in A}(a-b),$ and therefore $\sum_{a\in A} a=(m-1)b.$ Since \(p\nmid(m-1)\), this determines \(b\) uniquely. Thus every point of \(C_{x,L}\setminus L\)
lies on the unique line \(M\) through \(x\) and $b$. Hence \(C_{x,L}\subseteq L\cup M\), as required. \qed

\subsubsection*{Group of the maximal cliques $C_{x,L}$}

As before, consider the situation of a Desarguesian net with a line $L$
and a point $x$ outside. Put $A = x^\perp \cap L$.

\begin{Proposition} \label{group}
Let $G_{x,L}$ be the group of invertible linear transformations stabilizing
$L$ and $C_{x,L}$. Then $G_{x,L}$ is transitive on $C_{x,L} \setminus L$.
\end{Proposition}
\Proof
This is clear from the proof of Theorem~\ref{netcliq}. Indeed, if we take
$L = \Ff_r$, then $G_{x,L}$ consists of linear maps $z \mapsto cz+d$
with $c,d \in \Ff_r$ and $c\neq 0$, such that $cx+d \in  C_{x,L} \setminus L$.
\qed

\begin{Proposition}
The group of linear automorphisms of the net that preserve $L$ and $C_{x,L}$
is faithful on $L$, and induces on $L$ the group $G$ of invertible linear maps
$g$ on $L$ preserving $A$ such that if $g$ has a unique fixed point $c$,
then $c \in A$.
\end{Proposition}
\Proof
Consider the linear map $z \mapsto az+b$ with $a,b \in \Ff_r$, $a \ne 0$.
It has a unique fixed point when $a \ne 1$, and then that fixed point is
$\frac{-b}{a-1}$.
This map sends $x$ to a point $ax+b$ collinear with $x$ if and only if
$\frac{-b}{a-1} \in A$. Thus $G$ is the group of linear automorphisms
of the net stabilizing $L$ and $C_{x,L}$, and in particular is a group.
\qed

\begin{Proposition}\label{fixedpoints}
Let $L$, $C_{x,L}$ and $G$ be as above, and let $C$ the set of fixed points
of maps $z \mapsto az+b$ with $a \ne 1$ in $G$. If $C \ne \emptyset$,
then $C$ is a single $G$-orbit.
\end{Proposition}
\Proof
See \cite{HughesPiper73}, Theorem 4.25. \qed

\subsubsection*{Sizes of the maximal cliques $C_{x,L}$}

Since we are considering $C_{x,L}$, there is at least one line, i.e.,
$m \ge 1$.
For $m = 1$ or $m = 2$ or $m = r-1$ or $m = r$ one has $|C_{x,L}| = r$.
%
%
%
%
%
%
For $m = r+1$ one has $|C_{x,L}| = r^2$.
The following proposition gives a necessary condition on the possible size of $C_{x,L}$ for $2 < m < r-1$.

\begin{Proposition} \label{CxLsize}
Let $2 < m < r-1$. One has $|C_{x,L}| = m-1+h p^f$
where $h$ is a positive integer such that $h \divides {\rm gcd}(r-1,m-2,p^f-1)$, and $f$ is a nonnegative integer such that $p^f \divides {\rm gcd}(r,m-1)$.
\end{Proposition}
\Proof
We have that the set $A=C_{x,L} \cap L$ has size $m-1$ and Proposition~\ref{group} implies that $|C_{x,L} \setminus L|$ is the size
of a single orbit of a subgroup $G$ of ${\rm AGL}(1,r)$, hence
of the form $h p^f$ for some integers $h,f$, with $h \divides (r-1)$
and $p^f \divides r$. Moreover, $G$ preserves $A$.

Identify $L$ with $\Ff_r$.
The group ${\rm AGL}(1,r)$ of linear transformations
$\{z \mapsto az+b \mid a,b \in \Ff_r, a \ne 0 \}$ (of size $r(r-1)$)
acts as a group of automorphisms on the net.
It has the two orbits $\Ff_r$ and $\Ff_{r^2} \setminus \Ff_{r^{\vphantom{2}}}$,
so no nonidentity element has a fixed point outside $L$. It follows that
the size of the orbit $Gx$ of $x$ equals $|G|$. We can write $G= T \rtimes H$, where $T$ is a normal subgroup of $G$ of order $p^f$ consisting of the translations in $G$ and $H = G/T$ is cyclic, isomorphic to a subgroup of $\Ff_r^*$,
so that $h \divides (r-1)$ if $h = |H|$. The orbits of $T$ all have size $p^f$ and thus $p^f \divides {\rm gcd}(r,m-1)$.
If $h=1$, then we are done. Next assume that $h>1$. Let $g \colon z \mapsto az+b$ be an element of $G$ such that $a$ has order $h$. Its orbits on $A$ all have size $h$, except for the singleton
$\smash{\{\frac{-b}{a-1}\}}$. It follows that $h \divides {\rm gcd}(r-1,m-2)$.
The subgroup $T$ is invariant under conjugation by $g$,
and the map $z \mapsto z+t$ is conjugated to $z \mapsto z+at$,
so $h \divides (p^f-1)$. \qed

The next proposition provides a sufficient condition: it is useful in constructing a net such that $C_{x,L}$ has a prescribed size.

\begin{Proposition} \label{suff}
Let $2<m<r-1$. Let $h \divides {\rm gcd}(r-1,m-2,p^f-1)$ and $p^f \divides {\rm gcd}(r,m-1)$.  Suppose that there is a subset $A$ of $\Ff_r$ of size $m-1$ together with a subgroup $G$ of ${\rm AGL}(1,r)$ with order $hp^f$ such that\\
(i) $A$ is invariant under $G$ with the property that for each nonidentity $T_{a,s}\in G$, we have $a\in A$, but \\
(ii) $A$ is not invariant under any larger such group. \\
Then for $L=\Ff_r$ and $x \notin \Ff_r$, in the net defined by $x^\perp \cap L=A$, we have $|C_{x,L}|=m-1+hp^f$.
\end{Proposition}

\Proof From the proof of Theorem~\ref{netcliq} and Proposition~\ref{group}, we have $C_{x,L} \setminus A=\{g(x): g \in G\}$ and thus $C_{x,L}$ has size $m-1+|G|=m-1+hp^f$.
\qed

Next we determine possible sizes of $C_{x,L}$. First we consider the case $p \nmid (m-1)$. In this case we have $f = 0$ in Proposition~\ref{CxLsize}, and we show there are no additional conditions.

\begin{Theorem}\label{pnmidm-1}
Let $2 < m < r-1$ and $p \notdivides (m-1)$.
There exists a net of degree $m$ with a maximal clique $C_{x,L}$ of size $m-1+h$
if and only if $h \divides {\rm gcd}(r-1,m-2)$.
\end{Theorem}
\Proof
By Proposition~\ref{CxLsize}, it suffices to construct a maximal clique $C_{x,L}$ of size $m-1+h$ for each $h \divides {\rm gcd}(r-1,m-2)$. 

If $h > 1$, let $H$ be the subgroup of order $h$ of $\Ff_r^*$. By Proposition~\ref{suff} with the group $G=\{z \mapsto hz: h \in H\}$, it suffices to construct a subset $A$ of $\Ff_r$ of size $m-1$ such that $0 \in A$ and $A\setminus\{0\}$ is a union of $H$-orbits but not of any larger subgroup. One can pick the $m-2$ points
in $\beta^i H$ for $i=0, \ldots, \frac{m-2}{h}-1$,
where $\beta$ is a primitive root in $\Ff_r$.

Finally, assume that $h = 1$. By Proposition~\ref{suff}, we have to pick a set $A = \{a_1,\ldots,a_{m-1}\}$,
not invariant for any $T_{a,s}$ with $a \in A$ and $s\neq 1$.
If $A$ is invariant for $T_{a,s}$, then $A-a$
is invariant for multiplication by $s$, so $s \sum_{i=1}^{m-1} (a_i-a) = \sum_{i=1}^{m-1} (a_i-a)$,
and since $s \ne 1$ it follows that $\sum_{i=1}^{m-1} a_i = (m-1)a$.
Since $p \notdivides (m-1)$, this determines $a$ uniquely, and
$a$ must occur among the $a_i$. 
Thus, it suffices to pick $m-1$ distinct nonzero elements $a_i$ that sum to zero,
if that is possible.

\smallskip
{\footnotesize
Is it possible, given $n$ with $2 \le n \le r-3$, to pick $n$ distinct
nonzero elements in $\Ff_r$ that sum to zero? Since $r \ne 2$, all nonzero
elements of $\Ff_r$ sum to zero, so we may assume that $n \le \frac12 (r-1)$.
Let $A$ be an arbitrary subset of $\Ff_r\setminus \{0\}$ of size $n-2$, with sum $s$.
Then we can take the last two elements to be $b$ and $-b-s$
provided $b$ is not in $A \cup ({-}s{-}A) \cup \{0,-s\}$ and $2b \ne -s$.
If $p \ne 2$ then at most $r-2$ field elements must be avoided
and there is a possible choice for $b$.
If $p = 2$, then we need $s \ne 0$, which can be arranged unless
$n \in \{2,r-3\}$.
\par}


There remains the case \(p=2\) and \(m-1=r-3\). Choose
$A=\Ff_r\setminus\{0,1,\lambda\},$
where \(\lambda\in\Ff_r^\times\setminus\{1\}\) is not of order \(3\). Such a choice is possible
since \(r\neq 4\). We claim that \(A\) is not invariant under any nontrivial \(T_{a,s}\) with
\(a\in A\). Indeed, if \(A\) were invariant, then so would be its complement $B=\{0,1,\lambda\}.$
Since \(a\notin B\), the action of \(T_{a,s}\) on \(B\) has no fixed point, and hence is a
\(3\)-cycle. Thus \(T_{a,s}\) has order \(3\), so \(s\) has order \(3\). Since \(0\in B\), the
set \(B\) must be the orbit of \(0\), namely $B=\{0,T_{a,s}(0),T_{a,s}^2(0)\}=\{0,as,as^2\}.$
Hence the ratio of the two nonzero elements of \(B\) has order \(3\), contradicting the choice of
\(\lambda\). 
\qed


\medskip
For the Paley graphs of order $r^2$ 
the cliques $C_{x,L}$ are conjecturally the second largest
(if they are not the largest).
Theorem~\ref{pnmidm-1} implies that in general net graphs with the same
parameters, larger non-line cliques can occur (since cliques remain
cliques when parallel classes are added to a net with
smaller $m$).
%
%
\begin{Corollary}
Let $2 < m < r-1$ and $p \notdivides (m-1)$.
If ${\rm gcd}(r-1,m-2) = 1$, then $C_{x,L} = \{x\} \cup (x^\perp \cap L)$.
If ${\rm gcd}(r-1,m-2) = 2$, and both $L$ and $x^\perp \cap L$ are invariant under
$z \mapsto -z$, and $0 \in x^\perp \cap L$, then
$C_{x,L} = \{x,-x\} \cup (x^\perp \cap L)$. \qed
\end{Corollary}

\vspace{-0.3cm}
This generalizes the result by Baker et al.~\cite{Baker-et-al96} for
the Paley graph $P(r^2)$ (that has $m = \frac{r+1}{2}$
and hence ${\rm gcd}(r-1,m-2) = 2$ if $r \equiv 3$ (mod 4)
and ${\rm gcd}(r-1,m-2) = 1$ otherwise).

Martin \& Yip \cite{MartinYip24} proved for Peisert graphs $P^*(r^2)$
with $r \equiv 3$ (mod 4) that the cliques $\{x\} \cup (x^\perp \cap L)$
are maximal.
\medskip

So far, we have dealt with the case \(p \nmid (m-1)\). If $p \divides (m-1)$, then the
generic case is as before, but there are cases where $A$ cannot be
chosen in general position.

\begin{Proposition} \label{restrictions}
Let $2 < m < r-1$ and $p \divides (m-1)$ and
$|C_{x,L}| = m-1+hp^f$ as in Proposition \ref{CxLsize}.
Then \\
(a) if $m-1 = p^f$ then $h+1$ is a power of $p$, and \\
(b) if $m-1 = (h+1)p^f$ then $h+1$ is not a power of $p$, and \\
(c) if $m-1 = r - 2p^f$ then $h = 2$.
\end{Proposition}
\Proof
We follow the same notations as in the proof of Proposition \ref{CxLsize}.
The set $A$ is a union of $G$-orbits. 
If $h = 1$ then all orbits have size $p^f$.
If $h > 1$ then a single orbit has size $p^f$,
and all other orbits have size $|G|=h p^f$.
By Proposition~\ref{fixedpoints}, the short orbit is precisely the set of fixed points
of the nonidentity elements of $G$ with a fixed point.
We may take $G$ to be generated by a map $z \mapsto cz$
of order $h$ together with the translation subgroup $T$ of order $p^f$.
Identifying $L$ with $\Ff_r$ we may take $0 \in A$.
Let $T = \{ z \mapsto z+b \mid b \in B\}$.

(a) If $|A|=p^f$, then $A$ is a single $G$-orbit.
If $h = 1$ then $p = 2$, since $G$ also contains the map $z \mapsto -z$
(that fixes $0 \in A$). If $h > 1$ then $A$ is invariant
under ${\rm AGL}(1,F)$ where $F = \Ff_p(c)$. But then $h = |F|-1$.

(b) If $|A|=(h+1)p^f$ then $A$ is the union of two $G$-orbits
$A_0$ and $A_1$, where $|A_0|=p^f$. Let $0 \in A_0$ and pick $y \in A_1$.
If $h+1$ is a power of $p$, then $c$ is primitive in the subfield
$F = \Ff_p(c)$ of $\Ff_r$ and hence $A$ consists of the points $my+b$
with $m \in F$ and $b \in B$. Thus $A$ is invariant under $z \mapsto z+y$.
This is a contradiction, since $T$ is the full translation group preserving $A$.

(c) If $|A| = r-2p^f$, then the condition $h \divides {\rm gcd}(r-1,m-2,p^f-1)$ implies that $h \mid 2$. Assume that $h\neq 2$,
then $h = 1$ and $L \setminus A$ is the union of two $G$-orbits of size $p^f$.
If $p = 2$ then $L \setminus A$ is invariant under a larger
translation group.
If $p > 2$, we can label $L$ so that $L \setminus A$
is also invariant under $z \mapsto -z$ (that fixes $0 \in A$).
This is a contradiction in both cases.
\qed

We now show that conversely, given $m,h,f$ satisfying the above requirements,
there is a Desarguesian net of degree $m$ with a clique $C_{x,L}$
of size $m-1+h p^f$.

\begin{Theorem}
Let $2<m<r-1$ and $p \mid (m-1)$. There exists a net of degree $m$
with a maximal clique $C_{x,L}$ of size $m-1+hp^f$ if and only if
$p^f \divides \gcd(r,m-1)$, $h \divides \gcd(r-1,m-2,p^f-1)$ and
conditions (a), (b), (c) in Proposition \ref{restrictions} hold.
\end{Theorem}
\Proof
It suffices to find a subset $A$ of $\Ff_r$ of size $m-1$ together with a subgroup $G$ of ${\rm AGL}(1,r)$ with order $hp^f$ satisfying the two conditions in Proposition~\ref{suff}. First note that since $h \divides {\rm gcd}(r-1,p^f-1)$ and $p^f \divides r$, the group ${\rm AGL}(1,r)$ has subgroups $G$ of order $h p^f$ with normal subgroup $T$ of order $p^f$. Any such $G$ will do, except when $m-1 = p^f$, where $G$ has to be chosen suitably.

\medskip
Suppose $m-1 = p^f$. By assumption (Prop.~\ref{restrictions}(a))
$h+1 = p^d$ for some $d$, and since $h \divides (p^f-1)$, we have
$d \divides f$. Let $A$ be a $\Ff_{p^d}$-subspace
of $\Ff_r$ of size $p^f$ that does not contain a subfield of $\Ff_r$
larger than $\Ff_{p^d}$. Now $G$ is ${\rm AGL}(1,p^f)$ and $A$ is
not invariant under a larger subgroup of ${\rm AGL}(1,r)$.

\smallskip
{\footnotesize
Given finite fields $K,L$ with $K < L$, say $|K| = q$, $|L| = q^m$,
can one pick a $K$-subspace of $L$ of given dimension
not containing any intermediate field?
It suffices to find such a hyperplane. There are
$q^{m-2}+\cdots+q+1$ hyperplanes on $K$
and fewer than $q^{m-i}$ on a subfield of size $q^i$.
So there are hyperplanes containing $K$ but no larger subfield.
\par}
\smallskip

Now let $m-1 \ne p^f$, and consider an arbitrary subgroup $G$
of order $h p^f$ of ${\rm AGL}(1,r)$. We show that it is possible
to choose $A$ of size $m-1$ to be invariant under $G$ but under
no larger group.

\smallskip
\textbf{Case 1: $h>1$.}

Then $G$ contains an element $z \mapsto az+b$
with $a$ of order $h$. By relabeling $L$ we can obtain $b = 0$ in this case.
Now $G$ is generated by the homothety $z \mapsto az$ with $a$ of order $h$,
and the normal translation subgroup $T$ of order $p^f$.
Define $B$ by $T = \{z \mapsto z+b \mid b \in B \}$.
An element $g$ of $G$ looks like $z \mapsto a^i z+b$ with $0 \le i < h$
and $b \in B$, so that $B$ is invariant under multiplication by $a$.
If $g$ has a unique fixpoint, then $a^i \ne 1$ and the fixpoint
is $-b/(a^i-1)$. But multiplication by $a^i-1$ is injective and preserves $B$,
so the fixpoints of $G$ are in $B$, and hence in the $A$ we will choose.
For $A$ we choose $B$ of size $p^f$ (the short orbit)
and $(m-1-p^f)/(hp^f)$ further orbits of $G$. Then part (i) in Proposition~\ref{suff} is satisfied.

If $m-1 = (h+1)p^f$ then $A$ is the union of $B$ and one further $G$-orbit.
If a larger group preserves $A$ and $A$ contains
all fixpoints of its nonidentity elements, then $A$ is a single orbit
under the larger group and has size a power of $p$,
which was excluded by assumption (Prop.~\ref{restrictions}(b)).

If $m-1 = r-hp^f$, then $L \setminus A$ is a single $G$-orbit,
and no larger $G$ is possible.

Now, by the divisibility assumptions on $p^f$ and $h$, we can write $m-1 = (th+1)p^f$ and $r = (sh+1)p^f$ for some integers $s$ and $t$ with $2 \le t \le s-2$.
Then the number of choices for $A$ is ${s \choose t}$, while
there are fewer choices for a larger group, so that $A$ can be chosen
so as not to admit a larger group.

\smallskip
{\footnotesize
Indeed, if a larger group preserves $A$ then it must have a larger $h$
or a larger $f$ (or both). The number of choices for $A$ is
in the former case not more than $\sum_d {s/d \choose t/d}$,
where $d$ runs through the primes dividing both $s$ and $t$,
and in the latter case not more than $\smash{s {[(s-1)/p] \choose [t/p]}}$,
altogether fewer than $(s+\frac12t) {[s/2] \choose [t/2]}$.
If $3 \le t \le s-4$, this is smaller than ${s \choose t}$ and we are done.
%
%
Instead of choosing $A \setminus B$, we can choose $L \setminus A$.
So, without loss of generality, $t \le \frac12s$. This leaves the cases
$t=2$, $s \ge 4$ and $(s,t)=(6,3)$ to be examined.
The latter case is ok, and for $t=2$ only $p=2$ needs to be examined.
If $s$ is even, equality can hold only when $2h+1$ is a power of 2,
impossible. If $s$ is odd, equality can hold only when $h+1 = 2$,
but we are in the case $h > 1$ here.
\par}
\smallskip

\medskip
\textbf{Case 2: $h=1$.}


If $f = 0$ (that is, $|G|=1$), we have to choose $m-1$ points
not invariant for a nontrivial group. If $p \ne 2$,
or $p = 2$ and $m \equiv 1$ (mod 4), then any choice
of $m-1$ points not summing to 0 suffices.
%
%
If $p = 2$ and $m \equiv 3$ (mod 4), then $m-1 \ne 2, r-2$
(by Prop.~\ref{restrictions}(b,c)), and we can pick a set of $m+1$
elements summing to 0 and not closed under $(x,y,z) \mapsto x+y+z$
and suitably delete two of these elements.
%
%

If $m-1 = 2p^f$ then $p \ne 2$ by assumption. If a larger group preserves $A$,
and $A$ contains all fixpoints of its nonidentity elements, then $A$ is
a single orbit under the larger group and has size a power of $p$,
impossible.

If $m-1 = r-p^f$, then $A$ is the complement of a single orbit of size
$hp^f = p^f$, and no larger group can act.

By assumption $m-1 \ne r-2p^f$. Finally, let $r = sp^f$ and let $m-1 = tp^f$ with $3 \le t \le s-3$.
Then $A$ can be chosen in ${s \choose t}$ ways,
while there are fewer choices for a larger group.

\smallskip
{\footnotesize
Indeed, if a larger group preserves $A$ then it must have a larger $h$
or a larger $f$ (or both). The number of choices for $A$ is
in the former case not more than
$s \sum_d {(s-1)/d \choose (t-1)/d}$,
where $d$ runs through the primes dividing $(p^f-1,s-1,t-1)$,
and in the latter case not more than
$\frac{s-1}{p-1} {[s/p] \choose [t/p]}$,
altogether fewer than $s(\frac12(t-1)+1){[s/2] \choose [t/2]}$.
If $3 \le t \le s-6$, this is smaller than ${s \choose t}$
and we are done. As before, we may assume $t \le \frac12s$.
This leaves nine cases with $s \le 10$, $t \le 5$ to be examined,
and since $s$ must be a power of $p$, none survives.
\par}
\smallskip
\qed

\subsubsection*{Sizes of the maximal cliques $C_{x,L}$ in a fixed net}

For a fixed net, and a fixed choice of $L$, all maximal cliques $C_{x,L}$
have the same size, since the group stabilizing $L$ is transitive on
the complement of $L$.
On the other hand, one can find cases where $C_{x,L}$ and $C_{x,M}$ have
different sizes.

For odd $r$, let $H = (\Ff_r^*)^2$ be the subgroup of $\Ff_r^*$ of index 2,
so that $|H| = \frac{r-1}{2}$.
Let $L,M$ be two distinct lines on 0, and $x,y$ points on $L,M$
(respectively) other than 0.
Then $C = \{0\} \cup xH \cup yH$ is a clique of maximum size $r$
in a net of degree $\frac{r+3}{2}$, and this is $C_{y,L}$,
the maximal clique containing $\{y\} \cup (y^\perp \cap L)$,
for that net. Let $N$ be the line on $0$ and $y-x$.
If $r > 9$ then $|C_{y,N}| < |C_{y,L}|$:
$|C_{y,N}| = \frac{r+5}{2}$ for $r \equiv 1$ (mod 4) and
$|C_{y,N}| = \frac{r+3}{2}$ for $r \equiv 3$ (mod 4).

\smallskip
{\footnotesize
Indeed, $N_0 = y^\perp \cap N =
\{0,w\} \cup \{\frac{h}{h-1}w \mid h \in H, h \ne 1\}$, where $w = y-x$.
Since $\smash{\frac{h}{h-1}+\frac{h^{-1}}{h^{-1}-1}} = 1$ the map
$\smash{T_{\frac12,-1}} \,:\, z \mapsto w-z$ preserves~$N_0$.
This forces $a = \frac12 w$.
If $r \equiv 3$ (mod 4), then $a \notin N_0$ and $C_{y,N} = \{y\} \cup N_0$.
If $r \equiv 1$ (mod 4), then $-1 \in H$, hence $\frac12 w \in N_0$.
One may check that for $r \ne 9$ no other $T_{a,s}$ preserves $N_0$,
and hence $C_{y,N} = \{y,w-y\} \cup N_0$. 
\par}

\section{Numerical data} \label{numerical}
\subsection{Maximal cliques in small Paley graphs}
{\footnotesize
The table below gives the sizes of the maximal cliques in the Paley graphs
$P(r^2)$ for $r \le 47$. Exponents are the number of nonequivalent orbits
of this size under the full group of the graph.

\smallskip
\addtolength{\tabcolsep}{-2pt}
\noindent
\begin{tabular}{@{}c|l@{}}
$r$ & sizes \\
\hline\\[-9pt]
3 & $3^1$ \\
5 & $3^1$, $5^1$ \\
7 & $5^1$, $7^1$ \\
9 & $5^3$, $9^1$ \\
11 & $7^3$, $11^1$ \\
13 & $5^{10}$, $7^4$, $13^1$ \\
%
17 & $5^3$, $7^{41}$, $9^9$, $17^1$ \\
19 & $7^{25}$, $8^7$, $9^{17}$, $11^4$, $19^1$ \\
23 & $7^{85}$, $8^{108}$, $9^{80}$, $10^7$, $11^9$, $13^4$, $23^1$ \\
25 & $7^{405}$, $8^{226}$, $9^{49}$, $13^2$, $25^1$\\
27 & $7^{27}$, $8^{411}$, $9^{142}$, $10^{50}$, $11^{12}$, $15^2$, $27^1$ \\
29 & $7^{410}$, $8^{1584}$, $9^{2104}$, $10^{148}$, $11^{46}$, $13^1$,
  $15^2$, $29^1$ \\
31 & $7^{60}$, $8^{2004}$, $9^{2734}$, $10^{933}$, $11^{199}$, $12^{26}$,
  $13^{46}$, $17^2$, $31^1$ \\
37 & $7^{103}$, $8^{2505}$, $9^{21556}$, $10^{14002}$, $11^{5712}$,
  $12^{219}$, $13^{222}$, $19^2$, $37^1$ \\
41 & $7^{168}$, $8^{7801}$, $9^{104495}$, $10^{62070}$, $11^{9583}$,
  $12^{149}$, $13^{128}$, $14^{19}$, $21^2$, $41^1$ \\
43 & $7^{15}$, $8^{1748}$, $9^{54700}$, $10^{109127}$, $11^{54759}$,
  $12^{9785}$, $13^{1490}$, $14^{156}$, $15^{87}$, $17^{20}$, $23^2$, $43^1$ \\
47 & $7^{12}$, $8^{1097}$, $9^{125545}$, $10^{434029}$, $11^{210725}$,
  $12^{28533}$, $13^{4904}$, $14^{628}$, $15^{230}$, $16^{27}$, $17^{50}$,
  $25^2$, $47^1\hspace{-5pt}$
\end{tabular}\addtolength{\tabcolsep}{2pt}

\smallskip
For $r \equiv 3$ (mod 4) this confirms the values
from Kiermaier \& Kurz \cite{KiermaierKurz09}.

\medskip\noindent
The table below gives the same information for the
smallest maximal cliques for $r \le 73$.

\smallskip\noindent
\addtolength{\tabcolsep}{-2pt}
\begin{tabular}{@{~}c|cccccccccccccccccccccccccccccc@{}}
$r$ & 3 & 5 & 7 & 9 & 11 & 13 & 17 & 19 & 23 & 25 & 27 & 29 \\
\hline\\[-9pt]
size & $3^1$ & $3^1$ & $5^1$ & $5^3$ & $7^3$ & $5^{10}$ & $5^3$ & $7^{25}$ & $7^{85}$ & $7^{405}$ & $7^{27}$ & $7^{410}$ \\
$r$ & 31 & 37 & 41 & 43 & 47 & 49 & 53 & 59 & 61 & 67 & 71 & 73 \\
\hline\\[-9pt]
size & $7^{60}$ & $7^{103}$ & $7^{168}$ & $7^{15}$ & $7^{12}$ & $7^2$ & $8^{455}$ & $8^{113}$ & $7^1$ & $8^9$ & $9^{119319}$ & $9^{187566}$
\end{tabular}
\addtolength{\tabcolsep}{2pt}
\par}

\subsection{Maximal cliques in the Taylor extension of small Paley graphs}
Given a strongly regular graph $\Gamma$ on $v$ vertices with $k = 2\mu$,
its Taylor extension $\Sigma$ is a distance-regular graph on $2(v+1)$ vertices
with intersection array $\{v,v-k-1,1;\,1,v-k-1,v\}$
(cf.~\cite[\S1.5]{BCN}, \cite[\S1.2.7]{BvM}),
an antipodal 2-cover of the complete graph $K_{v+1}$.

Recall that the Taylor extension \(\Sigma\) of \(\Gamma\) is an antipodal 2-cover of the
graph join \(\infty \vee \Gamma\) such that the neighborhood of each of the two vertices of \(\Sigma\) corresponding to \(\infty\) induces a copy of \(\Gamma\). Thus a clique in \(\Gamma\) gives rise to a clique in \(\Sigma\) by adjoining \(\infty\), and conversely every clique in \(\Sigma\) containing \(\infty\) restricts to a clique in \(\Gamma\). In particular, maximal cliques in \(\Sigma\) correspond to maximal cliques in \(\Gamma\), with their sizes increased by \(1\).

For the Paley graph of order $q$, the Taylor extension
is distance-transitive on $2(q+1)$ vertices, with automorphism group
$2 {\times} P\Sigma L_2(q)$ (cf.~\cite[p.\,228]{BCN}).
It follows that the maximal cliques in $\Sigma$ have sizes
that are 1 larger than those in $\Gamma$, while the number of orbits
is smaller. Below a table for $q = r^2$, $r \le 49$.

{\smallskip\footnotesize
We see that the extra automorphisms of $\Delta = \Gamma(0)$
are those flipping the edge $0\infty$, for $\Gamma = \Sigma(\infty)$.
\par}

\medskip\noindent
{\footnotesize
\addtolength{\tabcolsep}{-2pt}
\begin{tabular}{@{~}c|l@{}}
$r$ & sizes \\
\hline\\[-9pt]
3 & $4^1$ \\
5 & $4^1$, $6^1$ \\
7 & $6^1$, $8^1$ \\
9 & $6^2$, $10^1$ \\
11 & $8^2$, $12^1$ \\
13 & $6^6$, $8^2$, $14^1$ \\
17 & $6^2$, $8^{14}$, $10^4$, $18^1$ \\
19 & $8^8$, $9^2$, $10^5$, $12^3$, $20^1$ \\
23 & $8^{22}$, $9^{16}$, $10^{15}$, $11^1$, $12^4$, $14^2$, $24^1$ \\
25 & $8^{84}$, $9^{29}$, $10^{15}$, $14^1$, $26^1$ \\
27 & $8^6$, $9^{50}$, $10^{24}$, $11^8$, $12^6$, $16^1$, $28^1$ \\
29 & $8^{85}$, $9^{180}$, $10^{307}$, $11^{18}$, $12^{11}$, $14^1$, $16^1$, $30^1$ \\
31 & $8^{17}$, $9^{232}$, $10^{324}$, $11^{96}$, $12^{43}$, $13^3$, $14^{13}$, $18^1$, $32^1$ \\
37 & $8^{31}$, $9^{281}$, $10^{2471}$, $11^{1288}$, $12^{640}$, $13^{21}$, $14^{36}$, $20^1$, $38^1$ \\
41 & $8^{42}$, $9^{871}$, $10^{11298}$, $11^{5705}$, $12^{1003}$, $13^{17}$, $14^{29}$, $15^3$, $22^1$, $42^1$ \\
43 & $8^7$, $9^{196}$, $10^{5715}$, $11^{10050}$, $12^{4935}$, $13^{840}$, $14^{182}$, $15^{15}$, $16^{19}$, $18^5$, $24^1$, $44^1$ \\
47 & $8^5$, $9^{125}$, $10^{12980}$, $11^{39699}$, $12^{18351}$, $13^{2388}$, $14^{516}$, $15^{60}$, $16^{38}$, $17^3$, $18^{12}$, $26^1$, $48^1$ \\
49 & $8^1$, $9^{27}$, $10^{8069}$, $11^{27777}$, $12^{19587}$, $13^{3738}$, $14^{812}$, $15^{56}$, $16^{16}$, $17^2$, $18^5$, $26^1$, $50^1$
\end{tabular}
\addtolength{\tabcolsep}{2pt}
\par}

\subsection{Maximal cliques in small Peisert graphs}
Sizes of the maximal cliques in $P^*(r^2)$ for $3 \le r \le 43$.

\smallskip\noindent
{\footnotesize
\addtolength{\tabcolsep}{-2pt}
\begin{tabular}{@{~}c|l@{}}
$r$ & sizes \\
\hline\\[-9pt]
3 & $3^1$ \\
7 & $4^1$, $7^1$ \\
9 & $5^1$, $9^1$ \\
11 & $5^7$, $6^2$, $11^1$ \\
19 & $6^1$, $7^{69}$, $8^{40}$, $9^{27}$, $10^3$, $19^1$ \\
23 & $6^1$, $7^{222}$, $8^{442}$, $9^{186}$, $10^{22}$,
     $11^1$, $12^1$, $23^1$ \\
27 & $7^{205}$, $8^{809}$, $9^{273}$, $10^{16}$, $11^2$, $14^1$, $27^1$ \\
31 & $7^{157}$, $8^{6099}$, $9^{7998}$, $10^{1629}$, $11^{113}$,
     $12^{11}$, $13^{11}$, $16^1$, $31^1$ \\
43 & $7^2$, $8^{4495}$, $9^{121241}$, $10^{258708}$, $11^{121126}$,
     $12^{21011}$, $13^{2196}$, $14^{195}$, $15^{45}$, $16^{19}$,
     $17^8$, $22^1$, $43^1\hspace{-10pt}$ \\
47 & $7^6$, $8^{10614}$, ..., $12^{30760}$, $13^{2248}$, $14^{291}$, $15^{99}$, $16^{44}$, $17^4$, $19^5$, $24^1$, $47^1$ \\
49 & $7^2$, $8^{6350}$, ..., $12^{20222}$, $13^{962}$, $14^{100}$, $15^{21}$, $16^1$, $17^2$ \\
59 & $8^{141}$, ..., $15^{4968}$, $16^{1261}$, $17^{324}$, $18^{120}$, $19^{114}$, $20^{15}$, $21^{11}$, $23^6$, $30^1$, $59^1$
\end{tabular}
\addtolength{\tabcolsep}{2pt}

\medskip
If $r \equiv 3$ (mod 4) the maximum cliques in $P^*(r^2)$ have size $r$.
The cliques $C_{x,L}$ are maximal of size $(r+1)/2$.
The maximum cliques in $P^*(7^4)$ have size 17. There are two orbits.
The smallest maximal cliques there have size 7.
There are two orbits, one of which is the subfield $\Ff_7$.
\par}

\medskip
Sizes of the maximal cliques in the sporadic Peisert graph $P^{**}(23^2)$.

\smallskip\noindent
{\footnotesize
\addtolength{\tabcolsep}{-2pt}
\begin{tabular}{@{~}c|l@{}}
$r$ & sizes \\
\hline\\[-9pt]
23 & $6^3$, $7^{362}$, $8^{448}$, $9^{87}$, $10^2$, $11^1$, $12^1$, $23^1$
\end{tabular}
\addtolength{\tabcolsep}{2pt}
\par}

\section*{Acknowledgments}
The second author is supported by Natural Science Foundation of Hebei Province (A2023205045). The third author and the fourth author thank Hebei Normal University for hospitality on their visit with the second author in August 2023, where this project initiated. The fourth author also thanks Greg Martin and Venkata Raghu Tej Pantangi for helpful discussions. The authors are also grateful to anonymous referees for their valuable comments and corrections.


\begin{thebibliography}{99}

\bibitem{AsgarliYip22}
S. Asgarli \& C. H. Yip,
{\it Van Lint-MacWilliams' conjecture and maximum cliques
in Cayley graphs over finite fields},
J. Combin. Th. (A) {\bf 192} (2022) 105667.

\bibitem{Baker-et-al96}
R. D. Baker, G. L. Ebert, J. Hemmeter \& A. J. Woldar,
{\it Maximal cliques in the Paley graph of square order},
J. Statist. Plann. Inference {\bf 56} (1996) 33-38.

\bibitem{Ball03}
S.~Ball,
{\it The number of directions determined by a function over a finite field},
J. Combin. Th. (A) {\bf 104} (2003) 341--350.

\bibitem{Blokhuis84}
A.~Blokhuis,
{\it On subsets of GF$(q^2)$ with square differences},
Indag. Math. {\bf 46} (1984) 369--372.

\bibitem{Blokhuis-et-al99}
A. Blokhuis, S. Ball, A. E. Brouwer, L. Storme \& T. Szőnyi,
{\it On the number of slopes of the graph of a function
defined on a finite field},
J. Combin. Th. (A) {\bf 86} (1999) 187--196.

\bibitem{BCN}
A. E. Brouwer, A. M. Cohen \& A. Neumaier,
{\it Distance-regular graphs},
Springer, 1989.

\bibitem{BvM}
A. E. Brouwer \& H. Van Maldeghem,
{\it Strongly regular graphs},
Cambridge Univ. Press, 2022.

\bibitem{Carlitz60}
L. Carlitz,
{\it A theorem on permutations in a finite field},
Proc. Amer. Math. Soc.~{\bf 11} (1960) 456--459.

\bibitem{Goryainov-et-al18}
S. Goryainov, V. V. Kabanov, L. Shalaginov \& A. Valuzhenich,
{\it On eigenfunctions and maximal cliques of Paley graphs of square order},
Finite Fields Appl. {\bf 52} (2018) 361--369.

\bibitem{Goryainov-et-al22}
S. Goryainov, A. Masley \& L. Shalaginov,
{\it On a correspondence between maximal cliques in Paley graphs
of square order},
Discr. Math. {\bf 345} (2022) 112853.

\bibitem{GoryainovYip24}
S. Goryainov \& C. H. Yip,
{\it Extremal Peisert-type graphs without the strict-EKR property},
J. Combin. Th. (A) {\bf 206} (2024) 105887.

\bibitem{HughesPiper73}
D. R. Hughes \& F. C. Piper,
{\it Projective planes},
Springer, 1973.

\bibitem{KiermaierKurz09}
M. Kiermaier \& S. Kurz,
{\it Maximal integral point sets in affine planes over finite fields},
Discr. Math. {\bf 309} (2009) 4564--4575.

\bibitem{MartinYip24}
G. Martin \& C. H. Yip,
{\it Distribution of power residues over shifted subfields and
maximal cliques in generalized Paley graphs},
Proc. Amer. Math. Soc. {\bf 153} (2025) 109--124.

\bibitem{MuzychukKovacs05}
M.~Muzychuk \& I.~Kovács,
{\it A solution of a problem of A.~E.~Brouwer},
Des. Codes Cryptogr.~{\bf 34} (2005) 249--264.

\bibitem{Mullin09}
N.~Mullin,
{\it Self-complementary arc-transitive graphs and their imposters},
M.~Sc. Thesis, Univ. of Waterloo, 2009.

\bibitem{Peisert01}
W. Peisert,
{\it All self-complementary symmetric graphs},
J. Algebra {\bf 240} (2001) 209--229.

\bibitem{Sziklai99}
P.~Sziklai,
{\it On subsets of ${\rm GF}(q^2)$ with $d$th power differences},
Discr. Math. {\bf 208/209} (1999) 547--555.

\bibitem{Y26}
C. H. Yip,
{\it Van Lint--MacWilliams' conjecture and maximum cliques in Cayley graphs over finite fields, II.},
J. Combin. Th. (A) {\bf 221} (2026) 106162.

\end{thebibliography}
\end{document}